\documentclass[a4paper]{amsart}

\usepackage{amsmath}
\usepackage{amssymb}
\usepackage{float}

\def \bbN {{\mathbb N}}

\def \bbQ {{\mathbb Q}}

\begin{document}

\title[Numerical verification of the Collatz conjecture]{Numerical verification of the Collatz conjecture for billion digit random numbers}

\author{Andreas-Stephan Elsenhans}

\address{School of Mathematics and Statistics\\ University of Syndey NSW 2006\\ Australia}
\address{Institut f\"ur Mathematik\\ Universit\"at W\"urzburg\\ Emil-Fischer-Stra\ss e 30\\ D-97074 W\"urzburg\\ German
y}
\email{stephan.elsenhans@mathematik.uni-wuerzburg.de}
\urladdr{https://www.mathematik.uni-wuerzburg.de/institut/personal/elsenhans.html}


\maketitle

\begin{abstract}
The Collatz conjecture, also known as the $3n+1$ problem, is one of the most 
popular open problems in number theory.

In this note, an algorithm for the verification of the Collatz conjecture 
is presented that works on a standard PC for numbers with up to ten billion decimal places.
\end{abstract}

\section{Introduction}

The Collatz conjecture asserts that, for any positive integer \( n_0 \), the sequence 
defined by the recurrence relation
\[
n_{k+1} =
\begin{cases}
\frac{n_k}{2} & \text{if } n_k \text{ is even}, \\
\frac{3n_k + 1}{2} & \text{if } n_k \text{ is odd},
\end{cases}
\]
will always eventually reach the value 1, regardless of the initial value \( n_0 \). Despite its simple 
formulation, the conjecture has resisted all attempts of a general proof, though it has been verified in~\cite{Barina} 
for all integers up to $2^{68}$.

To check the conjecture for the numbers $1,2,3,\ldots,N$, one computes the sequence for all initial values $n_0$ in the
given range until a number smaller than $n_0$ is reached. As this number has been checked before, the conjecture
is verified for all numbers in the initial range.

At this point, a large number of optimisations are possible. First, starting with an even number $n$ will 
lead to $\frac{n}{2}$. As this is smaller than the initial value, only odd numbers have to be checked.
Next, an odd number $n \equiv 1 \pmod{4}$ will lead to $\frac{3n+1}{4}$ after two iterations. As this is 
smaller than $n$, for every $n > 1$, it remains to check the conjecture for all integers in the residue class
$(3 \bmod 4)$. 

The above can be extended. For example, only the residue classes 7, 15, 27, and 31 modulo 32 need to be checked.
In other words, only these residue classes do not lead to a smaller number within the first five steps.

Our algorithm will use the same approach. I.e., it will condense a large number of iterations 
into a single step. Using this, one can verify the Collatz conjecture for random numbers with 
up to ten billion decimal places. 

Finally, some statistics on the number of iterations to reach 1 is given. It indicates that this number is
about $\log(n_0) / \log(2 / \sqrt{3})$ and the distribution is similar to a normal distribution.

\section{Collatz Polynomials}

\subsection*{Definition}
The function
$$
T \colon \bbN \rightarrow \bbN, n \mapsto 
\begin{cases}
\frac{n}{2} & \text{if } n \text{ is even}, \\
\frac{3n + 1}{2} & \text{if } n \text{ is odd},
\end{cases}
$$
is called the {\it Collatz-$T$-function}.

\subsection*{Fact}
{\it 
For every $N := 2^k$, there exist unique linear polynomials $C_{k,r} \in \bbQ[X]$\, ($r = 0,\ldots,N-1$) such that
$$
C_{k,r}(n) = T^{(k)}(n) 
$$
for all $n \equiv r \pmod{N}$.
Here, $T^{(k)}$ denotes the $k$-th iteration of the Collatz-$T$-function.}

\subsection*{Definition}
The polynomials $C_{k,r}$ are called {\it Collatz polynomials}.

\subsection*{Fact}
{\it 
The Collatz polynomials satisfy the relation
$$
C_{k+l,r}(X) =  C_{k,C_{l, r \bmod 2^l}(r) \bmod 2^k } (C_{l, r \bmod 2^l}(X)) 
$$
for all positive integers $k,l$ and $r \in \{ 0,\ldots,2^{k+l}-1 \}$. }

\subsection*{Application}
For $k = 2$, the Collatz polynomials are
$$
C_{2,0} = \frac{ X }{4}, ~
C_{2,1} = \frac{3 X + 1}{4} ,~
C_{2,2} = \frac{3 X + 2}{4},~
C_{2,3} = \frac{9 X + 5}{4}\, .
$$
One can read off these polynomials once more that only the residue class $(3 \bmod 4)$ 
contains integers larger than 4 that do not lead to a smaller value after two interations.
Further, for $k = 5$, one has
\begin{align*}
C_{5,3} &= \frac{9 X + 5}{32},~
C_{5,7} = \frac{81 X + 73}{32},~
C_{5,11} = \frac{27 X + 23}{32}, \\
C_{5,15} &= \frac{81 X + 65}{32},~
C_{5,19} = \frac{27 X + 31}{32},~
C_{5,23} = \frac{27 X + 19}{32},\\
C_{5,27} &= \frac{81 X + 85}{32},~
C_{5,31} = \frac{243 X + 211}{32}\, .
\end{align*}
One can read off these polynomials that only the residue classes $7, 15, 27,$ and $31$ modulo $32$ 
do not lead to a smaller value within the first five interations.

One can extend this easily and determine the residue classes modulo $2^k$ that do not lead to a smaller value within the first $k$ iterations. 
The number of such classes are listed in the table below: 
\medskip

\begin{center}
\begin{tabular}{c|cc||cc|c}
Module & Number of remaining & & &  Module & Number of remaining\\
       & residue classes     & & &         & residue classes     \\
\hline
2     & 1    & & &                    $2^2$   & 1 \\
$2^{3}$ & 2  & & &                    $2^{4}$ & 3 \\
$2^{5}$ & 4  & & &                    $2^{6}$ & 8 \\ 
$2^{7}$ & 13 & & &                    $2^{8}$ & 19 \\ 
$2^{9}$ & 38 & & & 
$2^{10}$ & 64  \\  
\vdots & \vdots & & & \vdots & \vdots \\ 
$2^{25}$ & 573162 & & & 
$2^{26}$ & 1037374  \\ 
$2^{27}$ & 1762293 & & & 
$2^{28}$ & 3524586  \\ 
$2^{29}$ & 6385637 & & &
$2^{30}$ & 12771274
\end{tabular}
\end{center}
\medskip

This sequence is well known ({\tt https://oeis.org/A076227}).
Note that our computation gives the residue classes, not just their number.

\section{An algorithm for large numbers}
\subsection*{Introduction}
It is well known that, when working with large numbers, special multiplication algorithms~\cite{Strassen} have to be used.
These are nowadays available in various computer algebra systems and the potential of such methods has to be utilised here.
In this work, the computer algebra system {\tt magma}~\cite{BCP} was used.

\subsection*{Data structure and code}
As all Collatz polynomials are of the shape $\frac{aX + b}{2^c}$, the triple $\langle a,b,c \rangle$ was used to represent them. 
This leads to the following basic functions:

\begin{verbatim}
// Evaluate the polynomial pol (given as a triple) at the integer x
// assuming that the result is an integer
function MyEval(p,x)
 return ShiftRight(x * p[1] + p[2],p[3]);
end function;

// Compute p(q(X)) for polynomials in triple representation
function MySubs(p,q)
 return <p[1]*q[1], p[1]*q[2] + ShiftLeft(p[2],q[3]), p[3]+q[3]>;
end function;

// Direct construction of the Collatz polynomial C_{l,r}
function PolyDirect(r,l)
 pol := <1,0,0>;
 for j := 1 to l do
  if IsEven(r) then
   r := r div 2;            pol[3] := pol[3] + 1;
  else
   r := (3*r+1) div 2;      pol := MySubs(<3,1,1>,pol);
  end if;
 end for;
 return pol;
end function;
\end{verbatim}
Using this, one can compute the Collatz polynomials for large modules by a binary splitting method:
\begin{verbatim}
// Store the fist Collatz polynomials in a list for fast access 
ColList := [[PolyDirect(i,l) : i in [0..2^l-1]] : l in [1..8]];

// Fast recursive construction of C_{l,r}
function PolyFast(r,l)
 if l le 8 then return ColList[l][(r mod 2^l) + 1]; end if;
 t1 := l div 2;    t2 := l - t1;
 p1 := PolyFast(ModByPowerOf2(r,t1), t1);
 p2 := PolyFast(ModByPowerOf2(MyEval(p1,ModByPowerOf2(r, l)),t2),t2);
 return MySubs(p2,p1);
end function;
\end{verbatim}

Now, one can check the Collatz conjecture for large random numbers 
by doing $\lfloor \log_2(n_k) \rfloor /2$ iterations in one step. 
Larger steps would result in larger and more expensive Collatz polynomials, whereas
smaller steps would result in a larger number of costly final evaluations.
\begin{verbatim}
 len := 10^5;  n := Random(10^len); cnt := 0;
 time repeat
  steps := Max(1,Ilog2(n) div 2);   cnt := cnt + steps;
  pol := PolyFast(n,steps);         n := MyEval(pol,n);
 until n eq 1;
 len, cnt;
\end{verbatim}

\subsection*{Results}
The code above was applied to random numbers with up to $10^{10}$ decimal places. In all cases, the iteration ends up in 1.
The run time for one number and a statistic on the number of iterations to reach 1 are listed in the table below. 
All computations are performed on one core of an Intel i7-12700 processor running at 4.7 GHz.

\medskip
\begin{center}
\begin{tabular}{r|r|r|r|r}
Number of         & Number of & Average Number       & Standard  & Time\\
decimal places    & examples  & of steps to reach 1  & Deviation &  \\
\hline
          10\,000 & 10000 &           160\,085 &   1531 & 0.02s \\
         100\,000 & 10000 &        1\,600\,728 &   4838 & 0.25s \\ 
      1\,000\,000 & 10000 &       16\,007\,868 &  15177 & 2.72s \\
     10\,000\,000 & 10000 &      160\,077\,568 &  48567 &  31s \\ 
    100\,000\,000 &  4000 &   1\,600\,785\,302 & 152901 &  353s \\ 
 1\,000\,000\,000 &   400 &  16\,007\,826\,613 & 516305 & 1.13h \\
10\,000\,000\,000 &     1 & 160\,079\,821\,246 & $-\!\!\!-\!\!\!-$ & 14h
\end{tabular}
\end{center}
\medskip
The statistics above indicates that the average number of steps to reach 1 is about $\log(n_0) / \log(2/\sqrt{3})$.
Further, the observed skewness is between $-0.022$ and $0.054$ and the observed kurtosis is between $2.967$ and $3.075$. 
An application of the Kolmogorov–Smirnov test does not detect a significant difference between the observed distribution
and a normal distribution.

\subsection*{Complexity}
The table above suggests that the run time of the  algorithm presented is $O(N^{1 + \epsilon})$, for any $\epsilon > 0$. Here, $N$
denotes the bit length of the input. 
A proof of such a statement would require an upper bound for the total stopping time of the Collatz sequence. Assuming that 
such a bound is of the form $C \cdot N$, one can apply the standard techniques to analyse algorithms using fast 
multiplication~\cite{Gathen} to prove the complexity statement.

\subsection*{Conclusion}
The investigation shows that the Collatz conjecture can be efficiently verified 
for random numbers with ten billion decimal places. 
This does not prove the conjecture. But it does illustrate what is possible using modern computer algebra systems. The
interested reader may consult~\cite{RXB} and~\cite{R2} for other attempts to verify the Collatz conjecture for large numbers.


\begin{thebibliography}{99}

\bibitem{Barina}
Barina, D.: \emph{Convergence verification of the Collatz problem},
The Journal of Supercomputing, Volume {\bf 77}, Issue 3 , (2021), 2681 -- 2688

\bibitem{BCP}
W.~Bosma, J.~Cannon, C.~Playoust: {\it The Magma algebra system.}
I. The user language, J.~Symbolic Comput.~{\bf 24} (1997), 235 -- 265.

\bibitem{Gathen}
J. von~zur~Gathen and J. Gerhard, {\it Modern computer algebra}, Cambridge Univ. Press, New York, 1999

\bibitem{RXB}
W. Ren, S. Li, R. Xiao and W. Bi: \emph{Collatz Conjecture for 2$\,\,\widehat{~}$100000-1 
Is True - Algorithms for Verifying Extremely Large Numbers},
in: 2018 IEEE SmartWorld, Ubiquitous Intelligence \& Computing, Advanced \& Trusted Computing, Scalable Computing \& Communications, 
Cloud \& Big Data Computing, Internet of People and Smart City Innovation (SmartWorld/SCALCOM/UIC/ATC/CBDCom/IOP/SCI), 
Guangzhou, China (2018), pp. 411--416 

\bibitem{R2}
W. Ren: \emph{Collatz Computation Sequence for Sufficient Large Integers is Random}, Preprint: {\tt https://eprint.iacr.org/2023/648.pdf}

\bibitem{Strassen}
Sch\"onhage, A.; Strassen, V.:
\emph{Schnelle Multiplikation gro\ss er Zahlen.}
Computing (Arch. Elektron. Rechnen)   {\bf 7} (1971), 281 -- 292.


\end{thebibliography}
\end{document}